\documentclass[a4paper,12pt]{amsart}
\usepackage{amsmath}
\usepackage{amsmath}
\usepackage{amssymb}
\usepackage{amsthm}
\usepackage{graphicx,color,wrapfig}

\theoremstyle{plain}
\newtheorem{lemma}{Lemma}
\newtheorem{thm}{Theorem}
\newtheorem{Definition}[lemma]{Definition}

\newtheorem{cor}[lemma]{Corollary}
\newtheorem{rem}[lemma]{Remark}

\newenvironment{Proof}[1][Proof]{
    \par
    \topsep 6pt plus 6pt
    \trivlist
    \item[\hskip\labelsep\bfseries #1.]\ignorespaces}%
    {\qed\endtrivlist
}

\newcommand\mpar[1]{\marginpar{\tiny \color{red} #1}}
\newcommand{\R}{\mathbb{R}}

\newcommand\ggreen[1]{{\color{green} #1}}

\begin{document}

\author{Hayk Mikayelyan, Henrik Shahgholian}

\title[Hopf's  lemma]{Hopf's lemma
for a class of singular/degenerate PDE-s}






\thanks{The second author is partially supported by Swedish Research Council}

\address{Hayk Mikayelyan \\Department of Mathematical Sciences \\
Xi'an Jiaotong-Liverpool University\\
Ren Ai Lu 111\\
215123 Suzhou \\ Jiangsu Province \\ 
PR China}
\email{Hayk.Mikayelyan@xjtlu.edu.cn}

\address{Henrik Shahgholian\\
Department of Mathematics\\ 
Royal Institute of Technology (KTH) \\
100 44 Stockholm \\
Sweden}
\email{henriksh@math.kth.se}

\begin{abstract}
This paper concerns Hopf's boundary point lemma, in certain $C^{1,Dini}$-type  domains, 
for a class of singular/degenerate  PDE-s, including $p$-Laplacian.  Using  geometric properties of levels sets for harmonic functions in convex rings, we  construct sub-solutions to our equations that play the role of a barrier from below. By comparison principle we then conclude  Hopf's lemma.

\end{abstract}

\maketitle

\section{Introduction}
In this paper we  consider Hopf's  lemma, in certain $C^{1,Dini}$-type  domains,  for  the following type of operators
\begin{equation}\label{eq}
\Delta_H u:=\text{div}(H(|\nabla u|)\nabla u), 
\end{equation}
where $H(t)=t^{-1}h(t)$, $h(0)=0$ and $h(t)$ is a monotone increasing continuous function. We will call the weak solutions of (\ref{eq}) $H$-harmonic. The additional condition we impose on the Dini modulus of continuity $\epsilon (t)$ is
that the function $t\epsilon(t)$ is convex (more discussion in Section \ref{secdin}). 

Equation (\ref{eq}) is the Euler-Lagrangian  of the functional  
\begin{equation}\label{funk}
I(v)=\int_{D} F(|\nabla v|)dx,
\end{equation} 
where $F(t)=\int_0^t h(\tau)d\tau$, $F\in C^1([0,\infty))$, $F(0)=F'(0)=0$ and $F$ is strictly convex.
Here and in the sequel $D\subset \R^n$ ($n\geq 2$) is a  domain.

This type of operators arise in applications  
dealing with  flows where  the flow-rate is proportional to 
$$
H(|\nabla u|)\nabla u=h(|\nabla u|)\frac{\nabla u}{|\nabla u|}.
$$

\subsection{Conditions on $F$}
Let us now list all the assumptions we impose on $h$, that can formally be divided into three groups:
\begin{itemize}
\item {\bf Physical} (see equations (\ref{cond1})), 
\item {\bf Coercive} (see equation (\ref{coercond}), or (\ref{Delta2})), 
\item {\bf Technical} (see equation (\ref{condition^})).
\end{itemize}

The so-called physical conditions have been already presented 
in the introduction.

\vspace{4mm}

\noindent \hspace{13mm}  $\bullet$  $h(0)=0$: with vanishing gradient flow vanishes;
\begin{equation}\label{cond1}
\text{\noindent $\bullet$ monotonicity of $h(t)$: larger gradient $\Longrightarrow$ more flow;}
\end{equation}
\hspace{13mm} $\bullet$  $H$ depends only on $|\nabla u|$: isotropy with respect to the position 

\hspace{13mm}and direction.
\vspace{3mm}

The coercivity condition is needed to assure the existence.
It is well known that this problem  is, in general, ill-posed if $h$ is bounded. The best illustration is the minimal surface equation 
$$
\text{div}\left( \frac{\nabla u}{\sqrt{1+|\nabla u|^2}}\right)
$$
in the annulus $B_2\backslash B_1$, with boundary data $0$ on $\partial B_2$ and $M$ on $\partial B_1$. For large enough $M$ the catenoid cannot reach the level $M$ without leaving the domain $B_2\backslash B_1$ and the equation has no solution in $W^{1,1}$.

To avoid this we can 
impose the strong condition
\begin{equation}\label{coercond}
c  t^{p-1} <h(t)< C t^{p-1}
\end{equation} 
for some $p>1$, 
which makes the application of the direct methods of the calculus of variations in the Sobolev space $W^{1,p}(D)$, $p>1$,
possible (see \cite{D}). 
In the special
case $H(t)=t^{p-2}$ we obtain the $p$-Laplacian.

\vspace{3mm}

Alternatively we can impose a weaker coercivity condition and work in Sobolev-Orlicz spaces $W^{1,F}(D)$. 
Let us shortly introduce these spaces following \cite{RR} (see also \cite{RR1}). 
The Orlicz norm is defined as follows
$$
\|u\|_F=\min\left\{M\,\,\Big|\,\,\int_D F\left(\frac{|u|}{M}\right)dx\leq F(1)\right\}.
$$
This norm defines the Banach space $L^F(D)$.

If we now denote by $g$ the inverse function of $h$ and define the Legendre transform of $F$ by 
$$
F^*(t)=\int_0^t g(\tau)d\tau,
$$ 
then assuming $h(1)=1$ one can easily prove using Young's inequality the generalization of H\"older's inequality
$$
\int_D uvdx\leq \|u\|_F\|v\|_{F^*}.
$$

If both $F$ and its 
Legendre transform $F^*$ satisfy the so-called $\Delta_2$ condition, i.e.,  there exists $t_0>0$, $C_0>0$ such that
\begin{equation}\label{Delta2}
F(2t)\leq C_0 F(t)\,\,\,\text{and}\,\,\, F^*(2t)\leq C_0 F^*(t)\,\,\,\text{for}\,\,\,t>t_0,
\end{equation}
then 
$$
(L^F(D))^*=L^{F^*}(D)
$$
and, in particular, $L^F(D)$ is reflexive, since $(F^*)^*=F$ (Theorem 10, page 112, \cite{RR}). 

Now analogously we can define the Sobolev-Orlicz space $W^{1,F}(D)$ by the
Sobolev-Orlicz norm as sum of Orlicz norms of $u$ and $|\nabla u|$.
Under $\Delta_2$ condition $W^{1,F}(D)$ will be reflexive and we can apply the direct methods
of the calculus of variations.

The condition (\ref{Delta2}) is somewhat weaker since it requires  polynomial of $F(t)$ 
growth only for large $t$ and leaves more freedom for the behavior for small $t$. 

\vspace{3mm}

In our proof the function
$$
R(t)=\frac{F''(t)}{F'(t)}
$$
plays an important role and we need the following technical condition on $R$.
We assume that 
for every positive, monotone increasing, bounded function $c(s)$ in $\R^+$,
 $0<c<c(s)<C<\infty$,
there exist constants $\alpha>0$ and $\beta>0$, depending on 
function $R$, and constants $c$ and $C$, such that
\begin{equation}\label{condition^}
\int_t^T R(c(s)s)ds\geq \alpha \int_{t}^{T}R(\beta s)ds.
\end{equation}

\begin{rem}\label{remcondR}
Condition (\ref{condition^}) is satisfied for more or less any "reasonable" 
function $R=F''/F'$. For
monotone decreasing $R$ one can take $\alpha=1$, $\beta=C$ (this covers the case $F(t)=t^p$, $p>1$).
The authors think that it is easier to check the condition (\ref{condition^}) for a given function $F$,
than to try to introduce a broad class of functions satisfying it. 
\end{rem}

\begin{Definition}[$H$-potential]
For two convex domains $K_1 \Subset K_2$
we call the minimizer
$u$ of
\begin{equation}\label{Hfunk}
J(v)=\int_{K_2\backslash K_1}    F(|\nabla v|) dx
\end{equation}
in the class of functions $\{v\in W^{1,F}_0(K_2)|\,v\equiv
1\,\,\,\text{on}\,\,\, K_1\}$ an $H$-potential (see \cite{RR}).
\end{Definition}


\subsection{Liapunov-Dini boundary}\label{secdin}

In the case of harmonic functions ($F(t)=t^2$) K.O. Widman (\cite{W}) using 
the Green representation was able to prove a Hopf-type result for domains with 
Liapunov-Dini boundary (see below). Moreover,
the following estimates for the second derivatives of the solution $v$ of uniformly elliptic equation
with H\"older continuous coefficients 
have been proved as well (see equation (2.4.1) in \cite{W})
\begin{equation}\label{widD2}
|D^2 v(x)|\leq C_D \frac{\epsilon(\delta(x))}{\delta(x)},
\end{equation}
where $\epsilon(t)$ is a Dini modulus of continuity (see \eqref{dini-cond} below) and $\delta(x)$ is the distance 
of the point $x$ from the boundary.
It is also shown that $C^{1,Dini}$ regularity is necessary for Hopf lemma in axially symmetric 
domains (see Remark 1 in \cite{W}). 

Since there is no Green representation for $p$-harmonic functions it is not 
possible to repeat Widman's direct estimates of the function's growth. Our proof is based on
barrier construction and works for the general operator $\Delta_H$ under some regularity assumptions on the boundary.


Let us present the definition of Liapunov-Dini surface following
\cite{W}.

\begin{Definition}\label{dini-cond}
A modulus of continuity $\epsilon(r)\searrow 0$ as ${r\to 0}$ is called
Dini modulus of continuity if $\int t^{-1}\epsilon(t)dt<\infty$.
\end{Definition}

\begin{Definition}
A Liapunov-Dini surface $S$ is a closed, bounded
$(n-1)-$dimensional surface satisfying the following conditions:

\noindent (a) At every point of $S$ there is a uniquely defined tangent 
hyper-plane, and thus also a normal.

\noindent (b) There exits a Dini modulus of continuity $\epsilon(t)$ such
that if $\beta$ is the angle between two normals, and $r$ is the
distance between their foot points, then the inequality
$\beta\leq\epsilon(r)$ holds.

\noindent (c) There is a constant $\rho_S>0$ such that for any point $x\in S$,
any line parallel to the normal at $x$ meets $S\cap B_{\rho_s}(x)$ at
most once.
\end{Definition}

In simple words the above definitions says, that the surface $S$ is locally the graph 
of a $C^{1,Dini}$ function in a ball of fixed radius.

Since in general the function $t\epsilon(t)$ is not convex we introduce a sub-class of Dini modules of continuity as follows.

\begin{Definition}\label{dini-conv}
A Dini modulus of continuity $\epsilon(r)$ is called convex-Dini if 
the function $t\epsilon(t)$ is convex.
\end{Definition}

\begin{rem}
Note that domains with $C^{1,\alpha}$ boundary are convex-Dini. 
The introduction of such a sub-class is necessary because our proof relies on the 
construction of barriers in convex rings with $C^{1,Dini}$ boundary. It it in general not true that 
for any Dini modulus of continuity $\epsilon(t)$ there is another Dini modulus of continuity 
$\tilde{\epsilon}(t)\geq\epsilon(t)$ such that $t\tilde{\epsilon}(t)$ is convex.
\end{rem}


\begin{rem}\label{inoutDini}
Note that if the Dini modulus of continuity is convex-Dini then a domain 
$D$ with Liapunov-Dini boundary satisfies a
kind of inner (outer) convex $C^{1,Dini}$ condition in the following
sense: There exists a convex Liapunov-Dini domain $K$ such that
for any point $x_0\in \partial D$ there exists a translation and
rotation $K_{x_0}$ of the domain $K$ satisfying
$$
K_{x_0}\subset D, \,\,(K_{x_0}\subset\R^n\backslash
D)\,\,\,\text{and}\,\,\,\partial K_{x_0}\cap\partial D=\{x_0\}.
$$
Moreover, we can take 
$$
K=K_{r_D}=B_{r_D}((0,\dots,0,r_D))\cap\{x\, | \, x_n > 2|x'|\epsilon(|x'|) \},
$$ 
where $x=(x',x_n)$, $r_D<\rho_{\partial D}/2$
and $\epsilon$ is the convex-Dini modulus of continuity. 
Without loss of generality (if necessary by modifying $\epsilon$ at ``corners'' with 
$\partial B_{r_D}((0,\dots,0,r_D))$) we can assume that 
$K$ has smooth boundary. Let us also observe that for any $a>0$ by taking $r_D$ small enough we can have 
$$
B_{(1-a)r_D}((0,\dots,0,r_D))\Subset K.
$$
Let us assume that for $ r_D$ 
\begin{equation}\label{34inside} 
B_{\frac{3}{4}r_D}((0,\dots,0,r_D))\Subset K.
\end{equation}
\end{rem}


In the sequel we will use as barriers the $H$-potentials in the convex rings 
\begin{equation}\label{nuclearbarriers} 
K\backslash B_{r_D/2}((0,\dots,0,r_D))\,\,\,\,
\text{and}\,\,\,\, 
B_{3r_D}((0,\dots,0,-r_D))\backslash (-K),
\end{equation}
where $-K$ is obtained from $K$ by symmetry with respect to the origin. We will refer to convex
rings (\ref{nuclearbarriers}) as 
inner and out convex rings.




\section{The Main Result and its proof}

The main result of this paper is the following extension of the Kjel-Ove Widman's result to a wider class of operators (\ref{eq}), for which  $p$-Laplacian is a particular case, in Liapunov-Dini domains with convex-Dini
modulus of continuity. 
As we will see later (Remark \ref{outerDini}) our result yields the boundary Harnack principle 
for $H$-harmonic functions in domains with convex-Dini boundary.

\begin{thm}\label{hopf-mopf}
Assume $u$ is an $H$-harmonic function in the domain $D$. Further
assume $0\in \partial D$, $ \partial D$ satisfies the inner 
convex Dini condition at $0$ and 
$$
 u(x) > u(0)\,\,\,\text{for all}\,\,\,x\in D.
 $$ 
Then there exist positive constants $r_0$ and $c$ such
that
$$
\max_{B_r\cap D} u(x)-u(0)>cr,
$$
for $0<r<r_0$.
\end{thm}


\begin{Proof}

Without loss of generality we can assume that the outer normal of $\partial D$ a the origin is $(0,\dots,0,-1)$.
Since for the solutions of (\ref{eq}) we have the maximum principle (see Appendix II) 
we need to construct a barrier in the inner convex ring from (\ref{nuclearbarriers}).

Actually we will construct barriers in arbitrary convex ring $K_2\backslash K_1$, where 
$K_1\Subset K_2$ are two convex domains with Liapunov-Dini
boundary. 

From the Hopf lemma for harmonic functions (\cite{W})
we know that if $\Delta w=0$ in $K_2\backslash K_1$, with 
boundary values $w=0$ on $\partial K_2$ and $w=1$ on $\partial K_1$ 
then $\nabla w \not= 0$ on $ \partial (K_2\backslash K_1) $. 
Now we will prove the existence of a convex,
smooth, monotone increasing function $f:[0,1]\to[0,1]$, $f(0)=0$, $f(1)=1$, $f'(0)>0$, $f'(1)<\infty$ such that
\begin{equation}\label{fcondmain}
\Delta_H f(w)\geq 0
\end{equation}
in $K_2\backslash K_1$. This will mean that the function $f(w)$ is a subsolution
for $\Delta_H$, has non-vanishing gradient at any boundary point and thus can 
be used as a barrier. 

We start by computing 
\begin{equation*}
\Delta_H f(w)=H(|\nabla f(w)|)\Delta f(w) + 
\frac{H'(|\nabla f(w)|)}{|\nabla f(w)|}\Delta_\infty f(w),
\end{equation*}
where $\Delta_\infty u=\nabla u D^2 u \nabla u$ is the $\infty$-Laplace operator.
Using 
$$
\Delta f(w)=f'(w)\Delta w+f''(w)|\nabla w|^2=f''(w)|\nabla w|^2,
$$ 
and 
$$
\Delta_\infty f(w)=(f'(w))^3\Delta_\infty w+(f'(w))^2f''(w)|\nabla w|^4,
$$
we arrive at
\begin{multline}\label{vorihashiv}
\Delta_H f(w)=H(f'(w)|\nabla w|)f''(w)|\nabla w|^2 + \\
\frac{H'(f'(w)|\nabla w|)}{f'(w)|\nabla w|}
\Big[(f'(w))^3\Delta_\infty w+(f'(w))^2f''(w)|\nabla w|^4\Big].
\end{multline}
We thus  need to find a function $f $,  such that $f'(t)>0$ for $t\in[0,1]$ and 
$\Delta_H f(w) \geq 0$. To comply with the latter we need  (see \eqref{vorihashiv})
\begin{multline*}
f''(w)|\nabla w|^2\Big[H(f'(w)|\nabla w|)+
f'(w)|\nabla w|H'(f'(w)|\nabla w|)\Big]\geq\\
-\frac{H'(f'(w)|\nabla w|)}{f'(w)|\nabla w|}
(f'(w))^3\Delta_\infty w,
\end{multline*}
which after substitution $H'(t)=\frac{F''(t)}{t}-\frac{F'(t)}{t^2}$ and
$H(t)+tH'(t)=F''(t)$ simplifies to
\begin{equation}\label{addedplus}
f''(w)\geq \left( 
\frac{F'(f'(w)|\nabla w|)}{F''(f'(w)|\nabla w|)}-
f'(w)|\nabla w| \right) |\nabla w|^{-5}\Delta_\infty w.
\end{equation}
Let us note that 
\begin{equation}\label{harmoniccurv}
0=\Delta w=\partial_{\nu\nu} w - (n-1)\kappa \partial_\nu w,
\end{equation}
where $\nu$ is the unit vector in the direction of $\nabla w$ and
$\kappa$ is the mean curvature of the level set; here we have used that the level sets 
of a positive harmonic potential are smooth.  From this we conclude
$$
\Delta_\infty w =(\partial_\nu w)^2 \partial_{\nu\nu}w= (n-1)\kappa |\nabla w|^3.
$$

Since the level sets of a harmonic potential in a convex ring are convex (see \cite{L}), the mean curvature and thus  $\Delta_\infty w$ is positive near the boundary. This along with $f'(w)>0$ and (\ref{addedplus}) (which is yet to be proven) implies  that it is enough to find a function  $f$ such that
\begin{equation}\label{anhavf}
f''(w)R(f'(w)|\nabla w|)\geq |\nabla w|^{-5}\Delta_\infty w,
\end{equation}
where
$R(t)=
\frac{F''(t)}{F'(t)}$.

\vspace{2mm}

The Hopf lemma proved in \cite{W} for harmonic functions yields
\begin{equation}\label{nabla-w}
0<c<|\nabla w(x)|<C<\infty,
\end{equation}
where the constants $c$ and $C$ depend only on the convex ring. Using this we easily obtain
$$
c\min (w(x), 1-w(x)) \leq \delta(x) \leq C\min (w(x), 1-w(x))
$$
and together with \eqref{widD2}   
\begin{multline*}
|\nabla w|^{-5}|\Delta_\infty w| \leq c^{-3} |D^2 w| \leq  c^{-3} C_D \min
\frac{\epsilon(\delta(x))}{\delta(x)}\leq  \\
c^{-3} C_D \frac{\epsilon(C\min (w, 1-w))}{c\min (w, 1-w)}=:
 \zeta(w),
\end{multline*}
where $\zeta(t)\in L^1([0,1])$ depend only on the convex ring. 

In order to have (\ref{anhavf}) we need to construct  a function $f$ such that
\begin{equation}\label{anhavfzeta}
f''(w)R(f'(w)|\nabla w|)\geq \zeta(w).
\end{equation}
For any $x\in K_2\backslash K_1$ let us denote by $\ell_x$ the gradient flow line of $w$ 
which contains $x$. Let us parametrize the curve $\ell_x$ by $w\in[0,1]$. 
We can now integrate (\ref{anhavfzeta}) on any $\ell_x$ in parameter $w$
\begin{equation}\label{anhavfzetaint}
\int_{w_1}^{w_2} f''(w)R(f'(w)|\nabla w|) dw\geq\int_{w_1}^{w_2} \zeta(w)dw,
\end{equation}

Observe that since 
the level sets of $w$ are convex the 
function $|\nabla w|$ on $\ell_x$ as a function of $w$ are monotone increasing 
(see equation (\ref{harmoniccurv})), 
but on the other hand we know that it is bounded by (\ref{nabla-w}). 
Thus we can apply our technical condition (\ref{condition^})
\begin{equation*}
\int_{w_1}^{w_2} f''(w)R(f'(w)|\nabla w|) dw=\int_{f'(w_1)}^{f'(w_2)}R(c(s) s)ds 
\geq\alpha\int_{f'(w_1)}^{f'(w_2)}R(\beta s)ds, 
\end{equation*} 
where $s=f'(w)$ and the function $c(s)= |\nabla w|(s)> 0$ is a monotone function such that $c<c(s)<C$.

If we now construct a function $f$ such that

\begin{equation}\label{anhavfinal}
\alpha\int_{f'(w_1)}^{f'(w_2)}R(\beta s)ds\geq \int_{w_1}^{w_2} \zeta(w)dw
\end{equation} 
for all $0<w_1<w_2<1$, then for this function $f$ the inequality (\ref{anhavfzetaint}) will be satisfied 
for all gradient flow lines $\ell_x$ and thus the inequality (\ref{anhavfzeta}) will be satisfied everywhere 
in $K_2\backslash K_1$, and we would be done.

Since $F'(0)=0$ and $F'(\infty)=\infty$, the function $R(t)=\frac{F''(t)}{F'(t)}$ is not 
integrable near zero and at $+\infty$, due to 
$$
\int_t^T R(\tau)d\tau=\log\frac{F'(T)}{F'(t)}.
$$  
As $w_1\to 0$ and $w_2\to 1$ the right hand side of (\ref{anhavfinal}) remains bounded ($\zeta\in L^1(0,1)$) 
and we can write 
\begin{multline}
\alpha\int_{f'(w_1)}^{f'(w_2)}R(\beta s)ds=\\
\frac{\alpha}{\beta}\left(\log F'(\beta f'(w_2))- \log F'(\beta f'(w_1))\right)
=\int_{w_1}^{w_2} \zeta(w)dw .
\end{multline}
Now we can take $f'(0)=m>0$ and construct
\begin{equation}\label{explicitform}
f'(w)=\beta^{-1}g\left(
F'(\beta m)
e^{\frac{\beta}{\alpha}\int_{0}^w \zeta(\tau)d\tau }
\right),
\end{equation}
where $g$ is the inverse function of $h=F'$ on $\R^+$.
Thus we obtained that (\ref{fcondmain}) is satistied for $f(w)=\int_0^w f'(\tau)d\tau $, where $f'$ is given by
(\ref{explicitform}). 

By changing the parameter $m\in (0,\infty)$ we can construct $f$ such that $f(1)$ is any positive number.

The proof of the theorem now will easily follow from applying the barrier $f(w)$ in the inner convex ring (\ref{nuclearbarriers}) with the parameter $m$ to be chosen such that 
$$
f(1)=\min_{x\in B_{r_D/2}((0,\dots,0,r_D)) } u(x)>0.
$$

\end{Proof}

\begin{rem}\label{scaleno}
Observe that if $f(w)$ is a sub-solution of (\ref{eq}) then in general we cannot say anything about the function  
$\alpha f(w)$, and we should construct the appropriate sub-solution by changing the parameter $m$ in (\ref{explicitform}).
\end{rem}

\begin{rem}\label{outerDini}
If the boundary value of a non-negative $H$-harmonic function $u$ vanishes in a neighborhood of $y$ and the 
boundary $\partial D$ satisfies the outer convex $C^{1,Dini}$ condition at $y$, then we can apply 
the super-solution barrier $f(1)-f(w)$ in the outer convex ring (see (\ref{nuclearbarriers})), and obtain the Lipschitz bound
$$
u(x)
\leq CM
\text{dist}(x, K) , 
$$
for $x\in B_{r_D}(y)\cap D$, where $r$, $C$ depending only on $D$ and $M= \max_{B_{r_D}(y)} u $.
\end{rem}

\begin{rem}
One can make the condition (\ref{condition^}) even weaker: there exists a constant $\alpha>0$
and a 
monotone increasing continuous function 
$$
L:\R^+\to \R^+,\,\,\,L(0)=0,\,\,\,L(\infty)=\infty
$$ 
such that 
\begin{equation}\label{generalcond}
\int_t^T R(c(s)s)ds\geq \alpha \int_{L(t)}^{L(T)}R(s)ds.
\end{equation}
For functions $F$ satisfying (\ref{generalcond}) one will obtain
$$
f'(t)=L^{-1}g\left(
F'(L(m))
e^{\frac{1}{\alpha}\int_{0}^t \zeta(\tau)d\tau }
\right),
$$
where $g$ is the same as in (\ref{explicitform}).
\end{rem}

\vspace{5mm}



\section{Appendix I: A comparison principle}\label{App-I}

The comparison  principle for $p$-harmonic functions is well known (see \cite{HKM}), 
but for the general operator $\Delta_H$ we could not find a reference. Therefore we shall  present a proof of this.

\begin{thm}\label{compprinc}
Let $u$ be a weak solution of (\ref{eq})and $v$ be its weak sub-solution in the domain $D$ with $C^1$ boundary.
Further let $v\leq u$ on $\partial D$ in the sense of trace operator. Then $v\leq u$ in $D$.   
\end{thm}
\begin{Proof}
Let us denote by $F^*$ the Legandre transform of $F$.
Observe that $g(t)=(F^*(t))'$ is the inverse function of
the function $h(t)=F'(t)$. By Young's inequality
$$
ab\leq F(a)+F^*(b)
$$
and the equality holds if and only if $b=h(a)$.

If $u$ is weak solution of (\ref{eq}) then from the convexity of $F$ it follows that 
$$
\int_D F(|\nabla u|)dx\leq \int_D F(|\nabla w|)dx
$$
for any $w$ such that $u-w\in W_0^{1,p}(D)$. Otherwise 
\begin{multline*}
\int_D F(|\nabla (u+t(w-u))|)\leq (1-t)\int_D F(|\nabla u|)dx+ t\int_D F(|\nabla w|)dx<\\
\int_D F(|\nabla u|)dx- \epsilon t
\end{multline*}
for some $\epsilon>0$ and differentiating in $t$ we obtain
\begin{equation}\label{sub-sol}
\int_D H(|\nabla u|)\nabla u\nabla(u-w)dx<0,
\end{equation}
which gives a contradiction after approximating $u-w\in W_0^{1,p}(D)$ by a test function $\phi\in C_0^\infty(D)$.

Let us now assume that $v\nleq u$ and take as a test function $\psi=(v-u)^+$. By the definition of the sub-solution
$$
\int_D H(|\nabla v|)\nabla v \nabla \psi dx \leq 0.
$$
Thus 
$$
\int_{D_1} H(|\nabla v|)\nabla v \nabla v dx \leq 
\int_{D_1} H(|\nabla v|)\nabla v \nabla u dx ,
$$
where $D_1=\text{supp} \psi\subset D$. Using now Young's inequality we obtain
\begin{multline}\label{ineq-1}
\int_{D_1} H(|\nabla v|)\nabla v \nabla udx \leq
\int_{D_1} h(|\nabla v|) |\nabla u|dx \leq\\
\int_{D_1} F(|\nabla u|)dx+
\int_{D_1} F^*( h(|\nabla v|))dx.
\end{multline}
Since $h(t)=F'(t)=tH(t)$, and Young's inequality is an equality for $b=h(a)$, we deduce
\begin{multline}\label{ineq-2}
\int_{D_1} H(|\nabla v|)\nabla v \nabla vdx =
\int_{D_1} h(|\nabla v|) |\nabla v|dx =\\
\int_{D_1} F(|\nabla v|)dx+
\int_{D_1} F^*( h(|\nabla v|))dx.
\end{multline}
By \eqref{sub-sol}--\eqref{ineq-2} we arrive at 
$$
\int_{D_1} F(|\nabla v|)dx\leq \int_{D_1} F(|\nabla u|)dx,
$$
where the inequality is strict unless $\nabla u=\nabla v$ a.e. in $D_1$; a contradiction in since
$u-v\in W^{1,p}_0(D_1)$.
\end{Proof}


\subsection*{Acknowledgment}
The first author is grateful to Stephan Luckhaus and Juan Luis Vazquez for inspiring discussions.


\end{document}